\documentclass[11pt]{amsart}
\usepackage{amsmath, amsthm, amsfonts, amssymb} 
\usepackage{url}

\newtheorem{theo}{Theorem}[section]

\newtheorem*{theoremA}{Theorem A}
\newtheorem*{theoremB}{Theorem B}
\newtheorem*{theoremC}{Theorem C}

\newtheorem{lem}[theo]{Lemma}
\newtheorem{prop}[theo]{Proposition}
\newtheorem{claim}[theo]{Claim}
\theoremstyle{definition}

\newtheorem{df}[theo]{Definition}

\newtheorem*{assumption}{Assumption A}

\newcommand{\R}{\mathbf{R}}
\newcommand{\C}{\mathbf{C}}
\newcommand{\Z}{\mathbf{Z}}
\newcommand{\F}{\mathbf{F}}

\newcommand{\N}{\mathbf{N}}

\newcommand{\Id}{\operatorname{Id}}
\newcommand{\Ad}{\operatorname{Ad}}

\newcommand{\Tr}{\operatorname{Tr}}

\newcommand{\Aut}{\operatorname{Aut}}
\newcommand{\SL}{\operatorname{SL}}

\newcommand{\op}{\operatorname{op}}

\newcommand{\cb}{\operatorname{cb}}

\newcommand{\Lim}{\operatorname{Lim}}
\newcommand{\ssolid}{\operatorname{s-solid}}
\newcommand{\ctr}{\operatorname{ctr}}
\newcommand{\SO}{\operatorname{SO}}
\newcommand{\SU}{\operatorname{SU}}
\begin{document}

\title[Strongly solid group factors]{Strongly solid group factors which are not interpolated free group factors}

\begin{abstract}
We give examples of non-amenable ICC groups $\Gamma$ with the Haagerup property, weakly amenable with constant $\Lambda_{\cb}(\Gamma) = 1$, for which we show that the associated ${\rm II_1}$ factors $L(\Gamma)$ are strongly solid, i.e. the normalizer of any diffuse amenable subalgebra $P \subset L(\Gamma)$ generates an amenable von Neumann algebra. Nevertheless, for these examples of groups $\Gamma$, $L(\Gamma)$ is not isomorphic to any interpolated free group factor $L(\F_t)$, for $1 < t \leq \infty$.
\end{abstract}

\author{Cyril Houdayer }

\address{CNRS-ENS Lyon \\
UMPA UMR 5669 \\
69364 Lyon cedex 7 \\
France}

\email{cyril.houdayer@umpa.ens-lyon.fr}

\subjclass[2000]{Primary 46L10; 46L54. Secondary 22D25; 37A20}

\keywords{Free group factors; Deformation/rigidity; Intertwining techniques; Group von Neumann algebra}

\maketitle

\section{Introduction}

In their remarkable recent work \cite{{ozawapopa}, {ozawapopaII}}, Ozawa and Popa obtained striking structural results for {\it group/group measure space} von Neumann algebras. For instance, they showed that if $\Gamma = \F_n$ is a free group, $2 \leq n \leq \infty$, \cite{ozawapopa} or $\Gamma$ is a lattice in $\SL(2, \R)$ or $\SL(2, \C)$ \cite{ozawapopaII}, then the group von Neumann algebra $L(\Gamma)$ is {\it strongly solid}, i.e. for any diffuse amenable subalgebra $P \subset L(\Gamma)$, the normalizer $\mathcal{N}_{L(\Gamma)}(P)$ of $P$ inside $L(\Gamma)$ generates an amenable von Neumann algebra. This strengthened two well-known indecomposability results for free group factors: Voiculescu's result in \cite{voiculescu96}, showing that $L(\F_n)$ has no Cartan subalgebra, which in fact exhibited the first examples of factors with no Cartan decomposition; and Ozawa's result in \cite{ozawa2003}, showing that the commutant in $L(\F_n)$ of any diffuse subalgebra must be amenable ($L(\F_n)$ are {\it solid}) which itself strengthened the indecomposability of $L(\F_n)$ into tensor product of ${\rm II_1}$ factors ({\it primeness} for free group factors) in \cite{Ge}.  

Lattices in $\SL(2, \R)$ are {\it measure equivalent} to $\F_2$, so the ${\rm II_1}$ factors arising from these lattices may be isomorphic to an amplification of a free group factor (i.e., to an interpolated free group factor \cite{{dykema94}, {radulescu1994}}). Ozawa and Popa asked in their paper (\cite{ozawapopa}, page 18) whether {\it any} strongly solid ${\rm II_1}$ factor $N$ with the complete metric approximation property (i.e. $\Lambda_{\cb}(N) = 1$) follows isomorphic to some interpolated free group factor $L(\F_t)$, $1 < t \leq \infty$, a question recently emphasized by Popa in his talks (see \cite{popa07} for instance). We answer this question in the negative by providing examples of non-amenable ICC (infinite conjugacy classes) groups $\Gamma$, for which we show that the corresponding ${\rm II_1}$ factors $L(\Gamma)$ are strongly solid, but nevertheless they are never isomorphic to any $L(\F_t)$, for $1 < t \leq \infty$. The groups $\Gamma$ are obtained as follows.

\begin{assumption}
Let $m \geq 2$. Let $\Lambda$ be an infinite countable {\it amenable} group and $\Upsilon_1, \dots, \Upsilon_m$ be infinite {\it abelian} groups together with group homomorphisms $\sigma^i : \Lambda \to \Aut(\Upsilon_i)$, for $i = 1, \dots, m$, satisfying:
\begin{equation*}
\forall i \in \{1, \dots, m\}, \forall h \in \Upsilon_i \backslash\{e\}, \{ g \in \Lambda : \sigma^i_g (h) = h\} = \{e\}.
\end{equation*} 
Write $\Upsilon = \Upsilon_1 \ast \cdots \ast \Upsilon_m$ for the free product. We denote by $\sigma : \Lambda \to \Aut(\Upsilon)$ the group homomorphism where $\Lambda$ acts {\it diagonally} on $\Upsilon$: if $h = h_1 \cdots h_n$ is a non-trivial reduced word in $\Upsilon$ with $h_j \in \Upsilon_{i_j}$, $i_1 \neq \cdots \neq i_n$, we have
\begin{equation*}
\sigma_g (h) = \sigma^{i_1}_g (h_1) \cdots \sigma^{i_n}_g (h_n).
\end{equation*}
We denote by $\Gamma = \Upsilon \rtimes \Lambda$ the corresponding semi-direct product. Observe that $\Gamma$ may be regarded as the amalgamated free product  
\begin{equation*}
\Upsilon_1 \rtimes \Lambda \ast_\Lambda \dots \ast_\Lambda \Upsilon_m \rtimes \Lambda.
\end{equation*} 
\end{assumption}

Assumption A is satisfied for instance if $\Lambda = \Z$, $\Upsilon_i = \bigoplus_\Z \Z / 2\Z$ and $\Lambda$ acts on $\Upsilon_i$ by {\it shift}. In this case, $\Upsilon_i \rtimes \Lambda$ is the {\it wreath product} $(\Z / 2\Z) \wr \Z$ and we have
\begin{equation*}
\Gamma = (\Z / 2\Z) \wr \Z \ast_\Z \cdots \ast_\Z (\Z / 2\Z) \wr \Z.
\end{equation*}
(More generally, instead of $\Z$ we can take $\Lambda$ to be any torsion-free amenable group in the above example.) Any group $\Gamma$ satisfying Assumption A is non-amenable, ICC and weakly amenable with constant $1$. (We refer to Proposition \ref{properties} for the properties of $\Gamma$.) The main result of this paper is the following:

\begin{theoremA}\label{A}
Let $\Gamma$ be a countable group satisfying Assumption A. Then the non-amenable ${\rm II_1}$ factor $L(\Gamma)$ is strongly solid. Moreover, $L(\Gamma)$ has the Haagerup property, the complete metric approximation property and is not isomorphic to any interpolated free group factor $L(\F_t)$, for $1 < t \leq \infty$.
\end{theoremA}

The fact that $L(\Gamma)$ is not isomorphic to an interpolated free group factor follows from Jung's result \cite{jung} ($L(\Gamma)$ is {\textquotedblleft strongly $1$-bounded\textquotedblright}). Peterson and Thom pointed out in \cite{PT} that a stronger property than {\it strong solidity} might hold true for the free group factors, namely any diffuse amenable subalgebra $P \subset L(\F_n)$ should have a unique maximal amenable extension. It is clear that our example does not satisfy this stronger property (take $P = L(\Lambda) \subset L(\Gamma)$).

The proof of Theorem A, following a {\textquotedblleft deformation/rigidity\textquotedblright} strategy, is a combination of the ideas and techniques in \cite{{ipp}, {ozawapopa}, {popamal1}}. 

A group $\Gamma$ which satisfies Assumption A has a vanishing first $\ell^2$-Betti number, so it is unclear whether $\Gamma$ admits a {\it proper} cocycle $(b, \pi , \mathcal{K})$ where the unitary representation $(\pi, \mathcal K)$ can be taken weakly contained in the left regular representation $\ell^2(\Gamma)$. Consequently, we cannot use Peterson's deformations arising from cocycles \cite{peterson4} in order to prove the strong solidity of $L(\Gamma)$ (see Theorem B in \cite{ozawapopaII}). Instead, regarding $L(\Gamma)$ as an amalgamated free product over $L(\Lambda)$
\begin{equation*}
L(\Gamma) = L(\Upsilon_1 \rtimes \Lambda) \ast_{L(\Lambda)} \cdots \ast_{L(\Lambda)} L(\Upsilon_m \rtimes \Lambda),
\end{equation*}
we will use the {\it free malleable deformation} by automorphisms $(\alpha_t)$ defined in \cite{ipp}.

The proof then consists in two parts. Let $\Gamma$ be a group satisfying Assumption A, and write $M = L(\Gamma)$, $M_i = L(\Upsilon_i \rtimes \Lambda)$. First, we show that given any amenable subalgebra $P \subset M$ such that $P$ does not embed into $M_i$ inside $M$, the normalizer $\mathcal{N}_M(P)$ generates an amenable von Neumann algebra (see Theorem \ref{step}). For this, we will exploit the facts that the deformation $(\alpha_t)$ does not converge uniformly on the unit ball $(P)_1$ and that $P \subset M$ is {\it weakly compact}, and use the technology from \cite{ozawapopa}. So if $P \subset M$ is diffuse, amenable such that $\mathcal{N}_M(P)''$ is not amenable, $P$ must embed into some $M_i$ inside $M$. Exploiting Popa's intertwining techniques, we prove that $\mathcal{N}_M(P)''$ is {\textquotedblleft captured\textquotedblright} in $M_i$ (which is amenable by assumption) and finally get a contradiction.

We then investigate the class $\mathcal{C}_{\ssolid}$ of countable groups $G$ which are weakly amenable with constant $1$ and for which the group von Neumann algebra $L(G)$ is strongly solid. This class contains all amenable groups, the free groups $\F_n$ \cite{ozawapopa}, lattices in $\SL(2, \R)$ or $\SL(2, \C)$ \cite{ozawapopaII} and all the groups $\Gamma$ which satisfy Assumption A by Theorem A. This class is obviously stable under taking subgroups. Our second result is the following:

\begin{theoremB}
The class $\mathcal{C}_{\ssolid}$ is stable under taking free products.
\end{theoremB}

Theorem B provides other new examples of groups which belong to the class $\mathcal{C}_{\ssolid}$. For instance, let $\Gamma$ be a group satisfying Assumption A. Then for any $n \geq 1$, the iterated free product $\Gamma^{\ast n}$ belongs to $\mathcal{C}_{\ssolid}$. Moreover, since $\Gamma$ has a vanishing first $\ell^2$-Betti number, i.e. $\beta_1(\Gamma) = 0$, it follows from \cite{l2} that $\Gamma$ and $\Gamma^{\ast n}$, $n \geq 2$, are not measure equivalent. Also, by \cite{cout}, $\Gamma^{\ast n}$ is never measure equivalent to a free group. Also, since $L(\Gamma)$ is strongly $1$-bounded \cite{jung} and is embeddable into $R^\omega$ (see Proposition \ref{properties}), $L(\Gamma)$ is {\it freely indecomposable}, i.e. $L(\Gamma)$ is not isomorphic to any free product of diffuse finite von Neumann algebras. In particular, $L(\Gamma) \ncong L(\Gamma^{\ast n})$, for any $n \geq 2$.

In \cite{ozawapopa}, Ozawa and Popa gave the first examples of group measure space constructions where the Cartan subalgebra $L^\infty(X, \mu) \subset L^\infty(X, \mu) \rtimes G$ in the crossed product ${\rm II_1}$ factor is unique up to {\it unitary} conjugacy. These examples include all the free ergodic {\it profinite} p.m.p. actions $G \curvearrowright (X, \mu)$ on the standard probability space where the group $G$ is a lattice of a product of Lie groups such as $\SO(n, 1)$, $\SU(n, 1)$ for $n \geq 2$, $\SL(2, \R)$ and $\SL(2, \C)$ \cite{{ozawapopa}, {ozawapopaII}}. Recall that a p.m.p. profinite action $G \curvearrowright (X, \mu)$ has the property that $L^\infty(X, \mu)$ is a limit of an increasing sequence $(Q_n)$ of $G$-invariant finite dimensional subalgebras. As noticed in \cite{ioana4}, a countable group $G$ admits free ergodic profinite p.m.p. actions iff $G$ is {\it residually finite}, i.e. $G$ has a decreasing sequence $(G_n)$ of normal finite index subgroups such that $\bigcap_n G_n = \{e\}$. Using the same techniques as in the proof of Theorem A, we obtain new examples of groups $G$ acting in a profinite way on $(X, \mu)$ such that the ${\rm II_1}$ factor $L^\infty(X, \mu) \rtimes G$ has a unique Cartan decomposition, in the spirit of results in \cite{{ozawapopa}, {ozawapopaII}}:

\begin{theoremC}
Let $G_1, G_2$ be weakly amenable groups with constant $1$, $|G_1| \geq 2$, $|G_2| \geq 3$. Denote by $G = G_1 \ast G_2$ their free product. Then $L(G)$ has no Cartan subalgebra. 

Assume moreover that $G_1, G_2$ are residually finite, so that $G$ is residually finite as well. Then, for any free ergodic profinite (or merely compact) p.m.p. action $G \curvearrowright (X, \mu)$ on the standard probability space, $L^\infty(X, \mu) \subset L^\infty(X, \mu) \rtimes G$ is the unique Cartan subalgebra up to unitary conjugacy. 
\end{theoremC}

Note that the first part of Theorem C was known when each $L(G_i)$ embeds into $R^\omega$ \cite{{jung}, {jung2}}. When combined with Gaboriau's results \cite{l2}, Theorem C shows that any ${\rm II_1}$ factor $L^\infty(X, \mu) \rtimes G$, where $G = G_1 \ast G_2$, $|G_1| \geq 2$, $|G_2| \geq 3$, $\Lambda_{\cb}(G_i) = 1$, $\beta_1(G_i) < \infty$, arising from a free ergodic profinite action $G \curvearrowright (X, \mu)$, has trivial fundamental group. Also, if $H_1 \ast H_2 = H \curvearrowright (X, \mu)$ is another such action, with $\beta_1(G) \neq \beta_1(H)$, then $L^\infty(X, \mu) \rtimes G \ncong L^\infty(X, \mu) \rtimes H$.

In Section \ref{preliminaries}, we review the necessary background on the intertwining techniques and weakly compact actions. The key result (see Theorem \ref{step}) is proven in Section \ref{keyresult}. Relying on this result and exploiting the intertwining techniques, we then prove the main results of the paper.

{\bf Acknowledgements.} The author thanks Prof. Sorin Popa for carefully reading an earlier version of this paper and for his useful comments. He also thanks Jesse Peterson for the stimulating discussions regarding this work during his visit at Vanderbilt University.

\section{Preliminaries}\label{preliminaries}

\subsection{Intertwining techniques}

We first recall some notation. Let $P \subset M$ be an inclusion of finite von Neumann algebras. The {\it normalizer of} $P$ {\it inside} $M$ is defined as
\begin{equation*}
\mathcal{N}_M(P) := \left\{ u \in \mathcal{U}(M) : \Ad(u) P = P \right\},
\end{equation*}
where $\Ad(u) = u \cdot u^*$. The inclusion $P \subset M$ is said to be {\it regular} if $\mathcal{N}_M(P)'' = M$. The {\it quasi-normalizer of} $P$ {\it inside} $M$ is defined as
\begin{equation*}
\mathcal{QN}_M(P) := \left\{ a \in M : \exists b_1, \dots, b_n \in M, aP \subset \sum_i Pb_i, Pa \subset \sum_i b_iP \right\}.
\end{equation*}
The inclusion $P \subset M$ is said to be {\it quasi-regular} if $\mathcal{QN}_M(P)'' = M$. Moreover,
\begin{equation*}
P' \cap M \subset \mathcal{N}_M(P)'' \subset \mathcal{QN}_M(P)''.
\end{equation*}
Let $A, B$ be finite von Neumann algebras. An $A$-$B$ {\it bimodule} $H$ is a complex (separable) Hilbert space $H$ together with two {\it commuting} normal $\ast$-representations $\pi_A : A \to \mathbf{B}(H)$, $\pi_B : B^{\op} \to \mathbf{B}(H)$. We shall intuitively write $a \xi b = \pi_A(x)\pi_B(y^{\op}) \xi$, $\forall x \in A, \forall y \in B, \forall \xi \in H$. We say that $H_B$ is {\it finitely generated} as a right $B$-module if $H_B$ is of the form $p L^2(B)^{\oplus n}$ for some projection $p \in \mathbf{M}_n(\C) \otimes B$. If $A$ is a finite von Neumann algebra, we will denote by $\ctr_A$ the {\it center-valued} trace of $A$.

In \cite{{popamal1}, {popa2001}}, Popa introduced a powerful tool to prove the unitary conjugacy of two von Neumann subalgebras of a tracial von Neumann algebra $(M, \tau)$. We will make intensively use of this technique. If $A, B \subset (M, \tau)$ are (possibly non-unital) von Neumann subalgebras, denote by $1_A$ (resp. $1_B$) the unit of $A$ (resp. $B$).

\begin{theo}[Popa, \cite{{popamal1}, {popa2001}}]\label{intertwining1}
Let $(M, \tau)$ be a finite von Neumann algebra. Let $A, B \subset M$ be possibly non-unital von Neumann subalgebras. The following are equivalent:
\begin{enumerate}
\item There exist $n \geq 1$, a possibly non-unital $\ast$-homomorphism $\psi : A \to \mathbf{M}_n(\C) \otimes B$ and a non-zero partial isometry $v \in \mathbf{M}_{1, n}(\C) \otimes 1_AM1_B$ such that $x v = v \psi(x)$, for any $x \in A$.

\item The bimodule $\vphantom{}_AL^2(1_AM1_B)_B$ contains a non-zero sub-bimodule $\vphantom{}_AH_B$ which is finitely generated as a right $B$-module. 

\item There is no sequence of unitaries $(u_k)$ in $A$ such that 
\begin{equation*}
\lim_{k \to \infty} \|E_B(a^* u_k b)\|_2 = 0, \forall a, b \in 1_A M 1_B.
\end{equation*}
\end{enumerate}
\end{theo}
If one of the previous equivalent conditions is satisfied, we shall say that $A$ {\it embeds into} $B$ {\it inside} $M$ and denote $A \preceq_M B$. For simplicity, we shall write $M^n := \mathbf{M}_n(\C) \otimes M$.

We make the following technical observation that will be useful in the next sections. Assume $A \preceq_M B$. Then there exist $n \geq 1$, a projection $p \in B^n$, a unital $\ast$-homomorphism $\psi : A \to pB^np$ and a non-zero partial isometry $v \in \mathbf{M}_{1, n}(\C) \otimes 1_AM1_B$ such that $x v = v \psi(x)$, for any $x \in A$. Note that $v^*v \leq p$ and $v^*v \in \psi(A)' \cap pM^np$. We do not have any control on the position of $v^*v$ in general. Nevertheless, we may assume that $p$ equals the  support projection of $E_{B^n}(v^*v)$. Indeed, write $q$ for the support of $E_{B^n}(v^*v)$. Since $v^*v \in \psi(A)' \cap pM^np$, it follows that $E_{B^n}(v^*v) \in \psi(A)' \cap pB^np$ and thus $q \in \psi(A)' \cap pB^np$. For any $x \in A$, 
\begin{equation*}
xvq = v \psi(x) q = vq \psi(x).
\end{equation*}
Of course, $vq$ is not a partial isometry, but $vq \neq 0$, because $E_{B^n}((vq)^*vq) = q E_{B^n}(v^*v) q = E_{B^n}(v^*v)$. Write $vq = w |vq|$ for the polar decomposition of $vq$. Define the unital $\ast$-homomorphism $\theta : A \to qB^nq$ by $\theta(x) = \psi(x)q$, for any $x \in A$. It follows that $w$ is a non-zero partial isometry in $\mathbf{M}_{1, n}(\C) \otimes 1_AM1_B$ such that $xw = w\theta(x)$, for any $x \in A$. Finally the support projection of $E_{B^n}(w^*w)$ equals $q$.

If $Q \subset (M, \tau)$ is an inclusion of finite von Neumann algebras, we shall denote by $\langle M, e_Q\rangle$ the {\it basic construction} where $e_Q : L^2(M) \to L^2(Q)$ is the orthogonal projection which satisfies:
\begin{equation*}
e_Q x e_Q = E_Q(x) e_Q, \forall x \in M.
\end{equation*}
Note that $E_Q : M \to Q$ is the unique $\tau$-preserving faithful normal conditional expectation. The basic construction $\langle M, e_Q\rangle$ is a {\it semifinite} von Neumann algebra with semifinite faithful normal trace $\Tr$ defined by:
\begin{equation*}
\Tr(x e_Q y) = \tau(xy), \forall x, y \in M. 
\end{equation*}

\subsection{The complete metric approximation property}

\begin{df}[Haagerup, \cite{haa}]
A finite von Neumann algebra $(M, \tau)$ is said to have the {\it complete metric approximation property} (c.m.a.p.) if there exists a net $\Phi_n : M \to M$ of ($\tau$-preserving) normal finite rank completely bounded maps such that
\begin{enumerate}
\item $\lim_n \|\Phi_n(x) - x\|_2 = 0$, $\forall x \in M$;
\item $\lim_n \|\Phi_n\|_{\cb} = 1$.
\end{enumerate}
\end{df}
If $M$ has the c.m.a.p., then $pMp$ has the c.m.a.p., for any non-zero projection $p \in M$. It follows from \cite{CowlingHaagerup} that a countable group $\Gamma$ is weakly amenable with constant $\Lambda_{\cb}(\Gamma) = 1$ iff the group von Neumann algebra $L(\Gamma)$ has the c.m.a.p.

\begin{df}[Ozawa \& Popa, \cite{ozawapopa}]
Let $\Gamma$ be a discrete group, let $(P, \tau)$ be a finite von Neumann algebra and let $\sigma : \Gamma \curvearrowright P$ be a $\tau$-preserving action. The action is said to be {\it weakly compact} if there exists a net $(\eta_n)$ of unit vectors in $L^2(P \bar{\otimes} \bar{P})_+$ such that
\begin{enumerate}
\item $\lim_n \|\eta_n - (v \otimes \bar{v})\eta_n\|_2 = 0$, $\forall v \in \mathcal{U}(P)$; 
\item $\lim_n \|\eta_n - (\sigma_g \otimes \bar{\sigma}_g)\eta_n\|_2 = 0$, $\forall g \in \Gamma$;
\item $\langle (a \otimes 1)\eta_n, \eta_n \rangle = \tau(a) = \langle \eta_n, (1 \otimes \bar{a}) \eta_n \rangle$, $\forall a \in M, \forall n$.
\end{enumerate}
These conditions force $P$ to be amenable. A von Neumann algebra $P \subset M$ is said to be {\it weakly compact inside} $M$ if the action by conjugation $\mathcal{N}_M(P) \curvearrowright P$ is weakly compact.
\end{df}

\begin{theo}[Ozawa \& Popa, \cite{ozawapopa}]\label{weakcompact}
Let $M$ be a finite von Neumann algebra with the complete metric approximation property. Let $p \in M$ be a non-zero projection and let $P \subset pMp$ be an amenable von Neumann subalgebra. Then $P$ is weakly compact inside $pMp$.
\end{theo}

\section{An intermediate step}\label{keyresult}

\subsection{The malleable deformation for amalgamated free products}

We begin this section by recalling the free malleable deformation from \cite{ipp}. We fix some notation that we will be using throughout. For $i = 1, \dots, m$, let $(M_i, \tau_i)$ be a finite von Neumann algebra with a distinguished f.n. trace. Let $B \subset M_i$ be a common von Neumann subalgebra such that ${\tau_i}{|_B} = {\tau_j}{|_B}$, $\forall i, j \in \{1, \dots, m\}$. Write $M = M_1 \ast_B \cdots \ast_B M_m$ for the amalgamated free product over $B$. Set $\widetilde{M}_i = M_i \ast_B (B \bar{\otimes} L(\Z))$ and
\begin{eqnarray*}
\widetilde{M} & = & \widetilde{M}_1 \ast_B \cdots \ast_B \widetilde{M}_m \\
& = & M \ast_B (B \bar{\otimes} L(\F_m)).
\end{eqnarray*}
The trace on $\widetilde{M}$ will be denoted by $\tau$. In $\widetilde{M}_i$, denote by $u_i$ the Haar unitary generating $L(\Z)$. Let $f : \mathbf{T} \to ]-1, 1]$  be the  Borel function satisfying $\exp(\pi \sqrt{-1} f(z)) = z$, $\forall z \in \mathbf{T}$. Define $h_i = f(u_i)$ a selfadjoint element in $\widetilde{M}_i$ such that $\exp(\pi \sqrt{-1}h_i) = u_i$. Write $u_i^t = \exp(t \pi \sqrt{-1} h_i) \in \mathcal{U}(\widetilde{M}_i)$, for any $t \in \R$. Following \cite{ipp}, define the deformation $(\alpha_t)$ on $\widetilde{M} = \widetilde{M}_1 \ast_{B} \cdots \ast_{B} \widetilde{M}_m$ by:
\begin{equation*}
\alpha_t = (\Ad u_1^t) \ast_{B} \cdots \ast_{B} (\Ad u_m^t), \forall t \in \R.
\end{equation*}
Since $u_i^t \to 1$ strongly, as $t \to 0$, it is clear that $\alpha_t \to \Id$ pointwise in $\| \cdot \|_2$, as $t \to 0$. Moreover define the period-$2$ automorphism $\beta$ on $\widetilde{M} = M \ast_{B} (B \bar{\otimes} L(\F_m))$ by:
\begin{eqnarray*}
\beta(x) & = & x, \forall x \in M, \\
\beta(u_i) & = & u_i^*, \forall i \in \{1, \dots, m\}. 
\end{eqnarray*}
It was proven in \cite{ipp} that the deformation $(\alpha_t, \beta)$ is {\it s-malleable}:
\begin{equation*}
\alpha_t \beta = \beta \alpha_{-t}, \forall t \in \R.
\end{equation*}
Moreover, note that $\alpha_t$ and $\beta$ are $\tau$-preserving, and $\alpha_t, \beta$ are equal to $\Id$ on $B$. We shall still denote by $\alpha_t$ and $\beta$ the corresponding isometries on $L^2(\widetilde{M})$. We remind at last that the s-malleable deformation $(\alpha_t, \beta)$ automatically features a certain {\it transversality} property.

\begin{prop}[Popa, \cite{popasup}]\label{transversality}
We keep the same notation as before. We have the following:
\begin{equation}\label{trans}
\|x - \alpha_{2t}(x)\|_2 \leq 2 \|\alpha_t(x) - (E_{M} \circ \alpha_t)(x)\|_2, \; \forall x \in M, \forall t > 0.
\end{equation}
\end{prop}

The following general result about intertwining subalgebras inside amalgamated free products will be a crucial tool in the next subsection (see also Theorem 4.2 in \cite{houdayer4} and Theorem 5.6 in \cite{houdayer3}).

\begin{theo}[Ioana, Peterson \& Popa, \cite{ipp}]\label{intertwiningipp}
Let $M = M_1 \ast_B \cdots \ast_B M_m$ be any amalgamated free product of finite von Neumann algebras $M_i$ over a common subalgebra $B$. Let $p \in M$ be a non-zero projection and $P \subset pMp$ be a von Neumann subalgebra. If the deformation $(\alpha_t)$ converges uniformly on the unit ball $(P)_1$, then there exists $i = 1, \dots, m$ such that $P \preceq_M M_i$.
\end{theo}

\subsection{The key result}

The following result will be a key ingredient in proving Theorems A, B, C.

\begin{theo}\label{step}
Let $M = M_1 \ast_B \cdots \ast_B M_m$ be an amalgamated free product with $B$ amenable. Let $Q_1, \dots, Q_k \subset M$ be amenable subalgebras such that the $M$-$M$ bimodule $L^2(\widetilde{M}) \ominus L^2(M)$ is a sub-bimodule of a multiple of $\bigoplus_{j = 1}^k L^2 \langle M, e_{Q_j}\rangle$. Let $p \in B$ be a non-zero projection and let $P \subset pMp$ be an amenable von Neumann subalgebra such that $P \npreceq_M M_i$, for $i = 1, \dots, m$. Let $\mathcal{G} \subset \mathcal{N}_{pMp}(P)$ be a subgroup such that the action by conjugation $\mathcal{G} \curvearrowright P$ is weakly compact. Then $\mathcal{G}''$ is amenable, thus AFD.
\end{theo}

\begin{proof}
The proof is similar to the ones of Theorems $4.9$ and $4.10$ in \cite{ozawapopa}. The only main difference is the fact that the deformation $(\alpha_t)$ is not compact over $B$. Instead we will use Theorem \ref{intertwiningipp}. Moreover, we need to pay attention to the fact that $P \subset M$ is {\it a priori} non-unital. We will nevertheless give a detailed proof for the sake of completeness. The symbol {\textquotedblleft Lim\textquotedblright} will be used for a state on $\ell^\infty(\N)$ which extends the ordinary limit.

Let $M = M_1 \ast_B \cdots \ast_B M_m$ be an amalgamated free product with $B$ amenable. Let $p \in B$ be a non-zero projection and let $P \subset pMp$ be an amenable von Neumann subalgebra such that $P \npreceq_M M_i$, for $i = 1, \dots, m$. Let $\mathcal{G} \subset \mathcal{N}_{pMp}(P)$ be a subgroup such that the action by conjugation $\mathcal{G} \curvearrowright P$ is weakly compact. We may and will assume that $\mathcal{U}(P) \subset \mathcal{G}$. Then there exists a net $(\eta_n)$ of vectors in $L^2(P \bar{\otimes} \bar{P})_+$ such that
\begin{enumerate}
\item $\lim_n \|\eta_n - (v \otimes \bar{v})\eta_n\|_2 = 0$, $\forall v \in \mathcal{U}(P)$; 
\item $\lim_n \|\eta_n - \Ad(u \otimes \bar{u})\eta_n\|_2 = 0$, $\forall u \in \mathcal{G}$;
\item $\langle (a \otimes 1)\eta_n, \eta_n\rangle = \tau(a) = \langle \eta_n, (1 \otimes \bar{a}) \eta_n \rangle$, $\forall a \in pMp, \forall n$.
\end{enumerate}
Note that $p$ is the unit of $P$, $\alpha_t(p) = p$, $\forall t \in \R$ and $\langle \cdot, \cdot\rangle$ denotes the inner product in $L^2(M \bar{\otimes} \bar{M})$. We regard $\eta_n \in L^2(M \bar{\otimes} \bar{M})_+$, and note that $(p \otimes 1) \eta_n (p \otimes 1) = \eta_n$, $(J \otimes \bar{J}) \eta_n = \eta_n$, where $J$ denotes the canonical anti-unitary on $L^2(M)$.

Fix $\varepsilon > 0$, $F \subset \mathcal{G}$ a finite subset and $z \in \mathcal{Z}(\mathcal{G}' \cap pMp)$ a non-zero projection. Observe that $z \in P' \cap pMp$. In particular, it follows that $Pz \npreceq_M M_i$, for $i = 1, \dots, m$. Using Theorem \ref{intertwiningipp}, we obtain that the deformation $(\alpha_t)$ does not converge uniformly on $(Pz)_1$. Since any selfadjoint element $x \in (Pz)_1$ can be written
\begin{equation*}
x = \frac12 \|x\|_\infty (u + u^*)
\end{equation*}
where $u \in \mathcal{U}(Pz)$, it follows that $(\alpha_t)$ does not converge uniformly on $\mathcal{U}(Pz)$ either. Combining this with the inequality $(\ref{trans})$ in Proposition \ref{transversality}, we get that there exist $0 < c < 1$, a sequence of positive reals $(t_k)$ and a sequence of unitaries $(u_k)$ in $\mathcal{U}(P)$ such that $\lim_{k \to \infty} t_k = 0$ and $\| \alpha_{t_k}(u_k z) - (E_M \circ \alpha_{t_k})(u_k z) \|_2 \geq c\| z\|_2$, $\forall k \in \N$. Since $\|\alpha_{t_k}(u_k z)\|_2 = \| z\|_2$, by Pythagora's theorem we obtain
\begin{equation}\label{key}
\|(E_M \circ \alpha_{t_k})(u_k z)\|_2 \leq \sqrt{1 - c^2}\| z\|_2, \forall k \in \N.
\end{equation}
Take $0 < \delta < \frac{1 - \sqrt{1 - c^2}}{6} \|z\|_2$. Choose and fix $k \in \N$ such that $\alpha = \alpha_{t_k}$ satisfies
\begin{eqnarray}\label{delta}
\|z - \alpha(z)\|_2 & \leq & \delta \\ 
\|u - \alpha(u)\|_2 & \leq & \varepsilon/6, \forall u \in F. \label{epsilon}
\end{eqnarray}
We set $v = u_k$. Define
\begin{eqnarray*}
\widetilde{\eta}_n & = & (\alpha \otimes 1)(\eta_n) \in L^2(\widetilde{M}) \bar{\otimes} L^2(\bar{M}) \\
\zeta_n & = & (e_M \otimes 1) (\widetilde{\eta}_n) \in L^2(M) \bar{\otimes} L^2(\bar{M}) \\
\zeta^\perp_n & = & \widetilde{\eta}_n - \zeta_n \in (L^2(\widetilde{M}) \ominus L^2(M)) \bar{\otimes} L^2(\bar{M}).
\end{eqnarray*}
Note that $(p \otimes 1) \zeta^\perp_n (p \otimes 1) = \zeta^\perp_n$ (since $\alpha(p) = p$), and moreover
\begin{equation}\label{norm2}
\|(xp \otimes 1) \widetilde{\eta}_n\|_2^2 = \tau(E_M(\alpha^{-1}(px^*xp))) = \|xp\|_2^2, \forall x \in \widetilde{M}.
\end{equation}
As in the proof of Theorem $4.9$ in \cite{ozawapopa}, $(\ref{epsilon}-\ref{norm2})$ yield for any $u \in F$,
\begin{eqnarray}\label{crucial}
\mathop{\Lim}_n \|[u \otimes \bar{u}, \zeta^\perp_n]\|_2 
& \leq & \mathop{\Lim}_n \|[u \otimes \bar{u}, \widetilde{\eta}_n]\|_2 \\ \nonumber
& \leq & \mathop{\Lim}_n \|(\alpha \otimes 1)([u \otimes \bar{u}, \eta_n])\|_2 + 2 \|u - \alpha(u) \|_2 \\
& < & \varepsilon/2. \nonumber
\end{eqnarray}
Moreover, $(\ref{key}-\ref{delta})$ and $(\ref{norm2})$ together with the choices of $t_k$ and $\delta$ yield 
\begin{equation}\label{crucial2}
\mathop{\Lim}_n \|(z \otimes 1) \zeta_n^\perp\|_2 > \delta.
\end{equation}
Indeed, suppose this is not the case. Noticing that $e_M z = z e_M$ (since $z \in M$) and $z v = v z$ (since $z \in \mathcal{Z}(\mathcal{G}' \cap pMp)$), with $v = u_k$, we have 
\begin{eqnarray*}
&& \mathop{\Lim}_n \|(z \otimes 1)\widetilde{\eta}_n - (e_M \alpha(v) z \otimes \bar{v})\zeta_n\|_2  \\
& \leq & \mathop{\Lim}_n \|(z \otimes 1)\widetilde{\eta}_n - (e_M \alpha(v) z \otimes \bar{v})\widetilde{\eta}_n\|_2 + \mathop{\Lim}_n \|(z \otimes 1) \zeta_n^\perp\|_2 \\
& \leq & \mathop{\Lim}_n \|(z \otimes 1)\widetilde{\eta}_n - (e_M z \alpha(v) \otimes \bar{v})\widetilde{\eta}_n\|_2 + \|[\alpha(v), z]\|_2 + \delta \\
& \leq & \mathop{\Lim}_n\|(z \otimes 1)\zeta_n^\perp\|_2 + \mathop{\Lim}_n \|\widetilde{\eta}_n - (\alpha(v) \otimes \bar{v})\widetilde{\eta}_n\|_2 + 2\|z - \alpha(z)\|_2 + \delta \\
& \leq & \mathop{\Lim}_n \|(\alpha \otimes 1)(\eta_n - (v \otimes \bar{v})\eta_n)\|_2 + 4\delta = 4 \delta.
\end{eqnarray*}
Thus, we get
\begin{eqnarray*}
\|(E_M \circ \alpha)(vz)\|_2 & \geq & \|(E_M \circ \alpha)(v)z\|_2 - \|z - \alpha(z)\|_2 \\
& \geq & \mathop{\Lim}_n \|((E_M \circ \alpha)(v)z \otimes \bar{v})\widetilde{\eta}_n\|_2 - \delta \\
& \geq & \mathop{\Lim}_n \|(e_M \otimes 1) ((E_M \circ \alpha)(v)z \otimes \bar{v})\widetilde{\eta}_n\|_2 - \delta \\
& = & \mathop{\Lim}_n \|(e_M \alpha(v)z \otimes \bar{v}) \zeta_n\|_2 -  \delta \\
& \geq & \mathop{\Lim}_n \|(z \otimes 1)\widetilde{\eta}_n\|_2 -  5\delta \\
& = & \|z\|_2 -  5\delta > \sqrt{1 - c^2} \|z\|_2,
\end{eqnarray*}
which is a contradiction according to $(\ref{key})$.

Choose $n$ large enough such that $\zeta = \zeta_n^\perp \in (L^2(\widetilde{M}) \ominus L^2(M)) \bar{\otimes} L^2(\bar{M})$ satisfies $(\ref{crucial})$ and $(\ref{crucial2})$, i.e. $\|[u \otimes \bar{u}, \zeta]\|_2 \leq \varepsilon/2$, $\forall u \in F$ and $\|(z \otimes 1)\zeta\|_2 \geq \delta$. Moreover for any $x \in M$, with $e_M^\perp = 1 - e_M$, we have
\begin{eqnarray*}
\| (xp \otimes 1) \zeta \|_2 & = & \|(xp \otimes 1) (e_M^\perp \otimes 1) \widetilde{\eta}_n\|_2 \\
& = & \|(e_M^\perp \otimes 1) (xp \otimes 1) \widetilde{\eta}_n\|_2 \\
& \leq & \| (xp \otimes 1) \widetilde{\eta}_n\|_2  = \|xp\|_2.
\end{eqnarray*}
By assumption, we may view $\zeta$ as a vector $(\zeta_i)$ in $\bigoplus_i L^2\langle M, e_{Q_{j(i)}}\rangle \bar{\otimes} L^2(\bar{M})$. Consider $\zeta_i \zeta_i^* \in L^1(\langle M, e_{Q_{j(i)}}\rangle \bar{\otimes} \bar{M})$, define $\xi = (\xi_i) \in \bigoplus_i L^2\langle M, e_{Q_{j(i)}}\rangle$, with $\xi_i = ( (\Id \otimes \tau)(\zeta_i\zeta_i^*) )^{1/2}$. We get $p \xi p = \xi$, $\|xp \xi\|_2 \leq \|xp\|_2$, $\forall x \in M$, $\|z \xi\|_2 = \|(z \otimes 1) \zeta\|_2 \geq \delta$, and $\|[u, \xi]\|_2 \leq \varepsilon$, $\forall u \in F$, as in the proof of Theorem 4.9 in \cite{ozawapopa}. (The last inequality follows from Powers-St\o rmer Inequality.)

Define now $\mathcal{G}_0 = \mathcal{G} + \mathbf{T}(1 - p) \subset \mathcal{U}(M)$ and $N = \mathcal{G}_0'' \subset M$. Hence $N$ is a {\it unital} von Neumann subalgebra of $M$. Observe that $N = \mathcal{G}'' + \C(1 - p)$ and $\mathcal{Z}(N' \cap M) = \mathcal{Z}( \mathcal{G}' \cap pMp) + \mathcal{Z}(M)(1 - p)$. Define $\xi' = \xi \oplus (1 - p)e_{B} \in (\bigoplus_i L^2\langle M, e_{Q_{j(i)}}\rangle) \oplus L^2\langle M, e_{B}\rangle$, $z' = z + z''(1 - p)$ where $z''$ is any projection in $\mathcal{Z}(M)$, and $F' = \{  u' = u + (1 - p) : u \in F\} \subset \mathcal{G}$. Note that $(1 - p)e_{B} = e_{B}(1 - p)$, since $p \in B$. For any $x \in M$,
\begin{eqnarray*}
\|x \xi'\|_2^2 & = & \|xp\xi \oplus x(1 - p)e_{B}\|_2^2 \\
& = & \|x p\xi\|_2^2 + \|x(1 - p)e_{B}\|_2^2 \\
& \leq & \|xp\|_2^2 + \|x(1 - p)\|_2^2 = \|x\|_2^2.
\end{eqnarray*}
Moreover, $\|z' \xi'\|_2 \geq \|z \xi\|_2 \geq \delta$ and $\|[u', \xi']\|_2 = \|[u, \xi]\|_2 \leq \varepsilon$, $\forall u' \in F'$. Applying Corollary 2.3 in \cite{ozawapopa}, we obtain projections $p_0, p_1, \dots, p_k \in \mathcal{Z}(N' \cap M)$ such that $p_0 + p_1 + \cdots + p_k = 1$ and $Np_0$ (resp. $Np_j$) is amenable relative to $B$ (resp. $Q_j$) inside $M$. Since $B, Q_1, \dots, Q_k$ are amenable, $Np_j$ is amenable, for any $j = 0, 1, \dots, k$. Since $N \subset Np_0 + Np_1 + \cdots + Np_k$ is a unital von Neumann subalgebra, $N$ is amenable and so is $\mathcal{G}''$.
\end{proof}

\section{$L(\Gamma)$ is strongly solid for $\Gamma$ satisfying Assumption A}\label{result}

\subsection{Properties of the group $\Gamma$ and its von Neumann algebra} Let $\Gamma$ be a countable group satisfying Assumption A. Recall that we may write $\Gamma$ in two different ways:
\begin{eqnarray*}
\Gamma & = & \Upsilon \rtimes \Lambda \\
& = & \Upsilon_1 \rtimes \Lambda \ast_\Lambda \cdots \ast_\Lambda \Upsilon_m \rtimes \Lambda,
\end{eqnarray*}
with $\Upsilon = \Upsilon_1 \ast \cdots \ast \Upsilon_m$. In the semi-direct product $\Upsilon \rtimes \Lambda$, the action $\sigma : \Lambda \to \Aut(\Upsilon)$ is implemented by {\it conjugation} with elements of $\Lambda$: $\sigma_g(h) = g h g^{-1}$, $\forall g \in \Lambda, \forall h \in \Upsilon$. We refer to \cite{{Gromov}, {Weiss}} for the notion of {\it sofic} groups.

\begin{prop}\label{properties}
Let $\Gamma$ be a group satisfying Assumption A. The following are true.
\begin{enumerate}
\item $\Gamma$ is non-amenable and ICC.
\item $\Gamma$ has the Haagerup property.
\item $\Gamma$ is weakly amenable with constant $\Lambda_{\cb}(\Gamma) = 1$.
\item $\Gamma$ is a sofic group. In particular, $L(\Gamma)$ is embeddable into $R^\omega$.
\item $\Gamma$ has a vanishing first $\ell^2$-Betti number, i.e. $\beta_1(\Gamma) = 0$.
\item The $\tau$-preserving action $\Lambda \curvearrowright L(\Upsilon)$ is mixing.
\item $L(\Gamma)$ is not isomorphic to any interpolated free group factor $L(\F_t)$, for $1 < t \leq \infty$. More generally, $L(\Gamma)$ is not isomorphic to any free product of diffuse finite von Neumann algebras.
\end{enumerate}
\end{prop}

\begin{proof}
Let $\Gamma$ be a group satisfying Assumption A. Since $\Upsilon = \Upsilon_1 \ast \cdots \ast \Upsilon_m$ is the free product, it follows that the group homomorphism $\Lambda \to \Aut(\Upsilon)$ satisfies
\begin{equation*}
\forall h \in \Upsilon \backslash\{e\}, \{ g \in \Lambda : ghg^{-1} = h\} = \{e\}
\end{equation*} 
as well. So if $h_1, h_2 \in \Upsilon \backslash \{e\}$, there exists at most one $g \in \Lambda$ such that $g h_1 g^{-1} = h_2$. Then, using this observation and the fact that $\Upsilon$ is a free product of infinite groups (in particular ICC), it is straightforward to check that $\Gamma$ is ICC as well. Thus, we get $(1)$. Using the above observation together with Kaplansky density theorem, $(6)$ follows easily. Being a free product of amenable groups, $\Upsilon$ has the Haagerup property \cite{jol}, is weakly amenable with Cowling-Haagerup constant $\Lambda_{\cb}(\Upsilon) = 1$ \cite{BP} and is sofic \cite{ES}. Since the sequence 
\begin{equation*}
1 \longrightarrow \Upsilon \longrightarrow \Gamma \longrightarrow \Lambda \longrightarrow 1
\end{equation*}
is exact and $\Lambda$ is amenable, $\Gamma$ has the Haagerup property, $\Lambda_{\cb}(\Gamma) = 1$ and $\Gamma$ is sofic as well (see \cite{{BO}, {haagerup}, {ES}}). We get $(2)$, $(3)$ and $(4)$. For $(5)$, regarding $\Gamma$ as an amalgamated free product of amenable groups over an infinite subgroup, it follows that $\Gamma$ has a vanishing first $\ell^2$-Betti number, i.e. $\beta_1(\Gamma) = 0$. In particular, $\Gamma$ cannot be embedded as a lattice in $\SL(2, \R)$ (by Gaboriau's result \cite{l2}). Observe that 
\begin{equation*}
L(\Gamma) = L(\Upsilon_1 \rtimes \Lambda) \ast_{L(\Lambda)} \cdots \ast_{L(\Lambda)} L(\Upsilon_m \rtimes \Lambda).
\end{equation*}
Since each $\Upsilon_i \rtimes \Lambda$ is amenable and ICC, the ${\rm II_1}$ factor $L(\Upsilon_i \rtimes \Lambda)$ is amenable hence AFD by Connes' result \cite{connes76}. Since $L(\Gamma)$ is an amalgamated free product of AFD ${\rm II_1}$ factors over a common diffuse subalgebra $L(\Lambda)$, it follows from Jung's result \cite{jung} that $L(\Gamma)$ is strongly $1$-bounded and therefore is never isomorphic to an interpolated free group factor $L(\F_t)$, for $1 < t \leq \infty$. More generally, since $L(\Gamma)$ is embeddable into $R^\omega$, $L(\Gamma)$ is not isomorphic to any free product of diffuse finite von Neumann algebras (see Lemma 3.7 in \cite{jung}).
\end{proof}

\subsection{The $M$-$M$ bimodule $L^2(\widetilde{M}) \ominus L^2(M)$}

Let $\Gamma = \Upsilon \rtimes \Lambda$ be a group satisfying Assumption A, with $\Upsilon = \Upsilon_1 \ast \cdots \ast \Upsilon_m$. Write $M_i = L(\Upsilon_i \rtimes \Lambda)$, $M = L(\Gamma)$ for the group von Neumann algebras, so that $M = M_1 \ast_{L(\Lambda)} \cdots \ast_{L(\Lambda)} M_m$. Denote by $\widetilde{\Gamma}$ the semi-direct product
\begin{equation*}
\widetilde{\Gamma} = (\Upsilon \ast \F_m) \rtimes \Lambda
\end{equation*}
where $\Lambda$ acts trivially on $\F_m$. Using the notation of Section \ref{keyresult}, it is clear that $\widetilde{M} = L(\widetilde{\Gamma})$.

\begin{prop}\label{bimodule}
Let $\Gamma$ be a countable group satisfying Assumption A. Let $M = L(\Gamma)$. Then as $M$-$M$ bimodules, we have
\begin{equation*}
L^2(\widetilde{M}) \ominus L^2(M) \cong \bigoplus (L^2\langle M, e_{L(\Lambda)}\rangle \oplus L^2\langle M, e_\C\rangle).
\end{equation*}
\end{prop}

\begin{proof}
Denote by 
\begin{enumerate}
\item $(u_g)_{g \in \Lambda}$ the canonical unitaries generating $L(\Lambda)$.
\item $(v_h)_{h \in \Upsilon}$ the ones generating $L(\Upsilon)$, where $\Upsilon = \Upsilon_1 \ast \cdots \ast \Upsilon_m$.
\item $(w_l)_{l \in \F_m}$ the ones generating $L(\F_m)$.
\end{enumerate}
Recall that if $x = \sum_{g \in \Lambda} x_g u_g \in M = L(\Upsilon) \rtimes \Lambda$, with $x_g \in L(\Upsilon)$, then
\begin{equation*}
E_{L(\Lambda)}(x) = \sum_{g \in \Lambda} \tau(x_g) u_g.
\end{equation*}
Denote by $\mathcal{I}$ the set of all the reduced words in $\Upsilon \ast \F_m$ of the form
\begin{equation*}
\xi = w_{l_1}v_{h_1} \cdots w_{l_k}v_{h_k}w_{l_{k + 1}},
\end{equation*}
where $k \geq 0$, $h_1, \dots, h_k, l_1, \dots, l_{k + 1} \neq e$. The {\it length} of $\xi$ is defined by $l(\xi) = 2k + 1$. For each $\xi \in \mathcal{I}$, denote by $\mathcal{H}_\xi$ the following $M$-$M$ bimodule:
\begin{equation*}
\mathcal{H}_\xi := \overline{M \xi M}.
\end{equation*}
By definition of the free product with amalgamation (see \cite{voiculescu92}), it follows that $\mathcal{H}_\xi$ is an $M$-$M$ sub-bimodule of $L^2(\widetilde{M}) \ominus L^2(M)$, for any $\xi \in \mathcal{I}$. Since $\Lambda$ normalizes $\Upsilon$ and $\F_m$ in $\widetilde{\Gamma}$, it follows that the $\mathcal{H}_\xi$'s generate $L^2(\widetilde{M}) \ominus L^2(M)$, i.e.
\begin{equation*}
L^2(\widetilde{M}) \ominus L^2(M) = \overline{\sum_{\xi \in \mathcal{I}} \mathcal{H}_\xi}.
\end{equation*}

For $\xi \in \mathcal{I}$, $\xi = w_{l_1}v_{h_1} \cdots w_{l_k}v_{h_k}w_{l_{k + 1}}$, write $g\cdot\xi := u_g \xi u_g^*$. This clearly defines an action of $\Lambda$ on the set $\mathcal{I}$, since
\begin{equation*}
g \cdot \xi = w_{l_1} v_{g h_1 g^{-1}} \cdots w_{l_k} v_{g h_k g^{-1}} w_{l_{k + 1}}.
\end{equation*}
Define now the {\it group stabilizer} of $\xi$ by 
\begin{equation*}
\Lambda_\xi := \{g \in \Lambda : g \cdot \xi = \xi\} = \bigcap_{j \leq k} \{g \in \Lambda : g h_j g^{-1} = h_j\}.
\end{equation*}
If $l(\xi) = 1$, then $\Lambda_\xi = \Lambda$, because $\Lambda$ commutes with $\F_m$. If $l(\xi) \geq 3$, then $\Lambda_\xi = \{e\}$, because $\{ g \in \Lambda : ghg^{-1} = h\} = \{e\}$, $\forall h \in \Upsilon \backslash\{e\}$. We will use the notation $\delta_{g, \Lambda_\xi}$ which equals $1$ if $g \in \Lambda_\xi$ and $0$ if $g \notin \Lambda_\xi$. We may identify $\Lambda \cdot \xi$ with the coset $\Lambda/\Lambda_\xi$. 

\begin{claim}\label{orthogonal}
Let $\xi, \eta \in \mathcal{I}$. If $\Lambda \cdot \xi = \Lambda \cdot \eta$, then $\mathcal{H}_\xi = \mathcal{H}_\eta$. If $\Lambda \cdot \xi \cap \Lambda \cdot \eta = \varnothing$, then $\mathcal{H}_\xi \perp_\tau \mathcal{H}_\eta$.
\end{claim}

\begin{proof}[Proof of Claim $\ref{orthogonal}$]
Let $\xi, \eta \in \mathcal{I}$. If $\Lambda \cdot \xi = \Lambda \cdot \eta$, then it is clear that $\mathcal{H}_\xi = \mathcal{H}_\eta$. Assume now that $\Lambda \cdot \xi \cap \Lambda \cdot \eta = \varnothing$. It suffices to show
\begin{equation*}
\langle x \xi y, z \eta t\rangle_\tau = 0,
\end{equation*}
for $x, y, z, t \in M$ with $x =  v_h u_\alpha$, $z = v_{h'} u_\gamma$, where $h, h' \in \Upsilon$, $\alpha, \gamma \in \Lambda$. We have
\begin{eqnarray*}
\langle x \xi y, z \eta t\rangle_\tau & = & \tau(t^*\eta^* z^* x \xi y) \\ 
& = & \tau(\eta^* z^* x \xi yt^*) \\ 
& = & \tau(\eta^* u_\gamma^* v_{h'}^* v_h u_\alpha \xi y t^*) \\
& = & \delta_{h, h'} \tau( \eta^* u_\gamma^* u_\alpha \xi y t^*) \\
& = & \delta_{h, h'} \tau(\eta^*(\gamma^{-1}\alpha \cdot \xi) u_\gamma^* u_\alpha y t^*) \\
& = & 0,
\end{eqnarray*}
because $u_\gamma^* u_\alpha y t^* \in M$ and the word $\eta^*(\gamma^{-1}\alpha \cdot \xi)$ is non-trivial, hence must contain at least one letter coming from $\F_m \backslash \{e\}$.
\end{proof}

\begin{claim}\label{isomorphism}
As $M$-$M$ bimodules, $\mathcal{H}_\xi \cong L^2\langle M, e_{L(\Lambda_\xi)}\rangle$, for any $\xi \in \mathcal{I}$.
\end{claim}

\begin{proof}[Proof of Claim $\ref{isomorphism}$]
We need to show that the map 
\begin{equation*}
\mathcal{H}_\xi \ni x \xi y \mapsto x e_{L(\Lambda_\xi)} y \in L^2\langle M, e_{L(\Lambda_\xi)}\rangle
\end{equation*}
preserves the inner products. It clearly suffices to show
\begin{equation}\label{inner}
\langle x \xi y, z \xi t\rangle_\tau = \langle x e_{L(\Lambda_\xi)} y, z e_{L(\Lambda_\xi)} t \rangle_{\Tr},
\end{equation}
for $x, y, z, t \in M$ with $x =  v_h u_\alpha$, $z = v_{h'} u_\gamma$, where $h, h' \in \Upsilon$, $\alpha, \gamma \in \Lambda$. On the left-hand side of $(\ref{inner})$, we have
\begin{eqnarray*}
\langle x \xi y, z \xi t\rangle_\tau 
& = &  \tau(t^*\xi^* z^* x \xi y) \\ 
& = & \tau(\xi^* z^* x \xi y t^*) \\
& = & \tau(\xi^* u_\gamma^* v_{h'}^* v_h u_\alpha \xi y t^*) \\
& = & \delta_{h, h'} \tau(\xi^* u_\gamma^* u_\alpha \xi y t^*) \\
& = & \delta_{h, h'} \tau(\xi^* (\gamma^{-1}\alpha \cdot \xi) u_\gamma^* u_\alpha y t^*) \\
& = & \delta_{\gamma^{-1}\alpha, \Lambda_\xi} \delta_{h, h'} \tau(u_\gamma^* u_\alpha y t^*),
\end{eqnarray*}
because if $\gamma^{-1}\alpha \notin \Lambda_\xi$, then the word $\xi^* (\gamma^{-1}\alpha \cdot \xi)$ is non-trivial and must contain at least one letter coming from $\F_m \backslash\{e\}$. On the right-hand side of $(\ref{inner})$, we have
\begin{eqnarray*}
\langle x e_{L(\Lambda_\xi)} y, z e_{L(\Lambda_\xi)} t \rangle_{\Tr} 
& = &  \Tr(t^* e_{L(\Lambda_\xi)} z^* x e_{L(\Lambda_\xi)} y) \\ 
& = & \Tr(t^* E_{L(\Lambda_\xi)}( z^* x) e_{L(\Lambda_\xi)} y) \\
& = & \tau(t^* E_{L(\Lambda_\xi)}( z^* x) y) \\
& = & \tau( E_{L(\Lambda_\xi)}( z^* x) y t^*) \\
& = & \tau( E_{L(\Lambda_\xi)}( (u_\gamma^* v_{h'}^* v_h u_\gamma) u_\gamma^* u_\alpha) y t^*) \\
& = & \tau(u_\gamma^* v_{h'}^* v_h u_\gamma) \tau( E_{L(\Lambda_\xi)}(u_\gamma^* u_\alpha) y t^*) \\
& = & \delta_{h, h'} \tau( E_{L(\Lambda_\xi)}(u_\gamma^* u_\alpha) y t^*) \\
& = & \delta_{\gamma^{-1}\alpha, \Lambda_\xi} \delta_{h, h'} \tau(u_\gamma^* u_\alpha y t^*).
\end{eqnarray*}
Consequently, $\langle x \xi y, z \xi t\rangle_\tau = \langle x e_{L(\Lambda_\xi)} y, z e_{L(\Lambda_\xi)} t \rangle_{\Tr}$.
\end{proof}
Using Claims \ref{orthogonal} and \ref{isomorphism} and the fact that for any $\xi \in \mathcal{I}$, $\Lambda_\xi = \{e\}$ or $\Lambda_\xi = \Lambda$, we are done.
\end{proof}

\subsection{Proof of Theorem A} Let's begin with a few easy observations first. Assume that $(N, \tau)$ is a finite von Neumann algebra with no amenable direct summand, i.e. $Nz$ is not amenable, $\forall z \in \mathcal{Z}(N)$, $z \neq 0$. Then for any non-zero projection $q \in N$, $qNq$ is non-amenable. Moreover, if $N$ has no amenable direct summand and $N \subset N_1$ is a unital inclusion of finite von Neumann algebras, then $N_1$ has no amenable direct summand either.

Let $\Gamma$ be a group satisfying Assumption A. Denote by $M = L(\Gamma)$ and by $M_i = L(\Upsilon_i \rtimes \Lambda)$ so that
\begin{equation*}
M = M_1 \ast_{L(\Lambda)} \cdots \ast_{L(\Lambda)} M_m.
\end{equation*}
Let $P \subset M$ be a diffuse amenable von Neumann subalgebra. By contradiction assume that $ \mathcal{N}_M(P)''$ is not amenable. Write $1 - z \in \mathcal{Z}(\mathcal{N}_M(P)'')$ for the maximal projection such that $\mathcal{N}_M(P)''(1 - z)$ is amenable. Then $z \neq 0$ and $\mathcal{N}_M(P)''z$ has no amenable direct summand. Notice that 
\begin{equation*}
\mathcal{N}_M(P)''z \subset \mathcal{N}_{zMz}(Pz)''.
\end{equation*}
Since this is a unital inclusion (with unit $z$), $\mathcal{N}_{zMz}(Pz)''$ has no amenable direct summand either. Since $L(\Lambda)$ is diffuse and $M$ is a ${\rm II_1}$ factor, there exist a projection $q \in L(\Lambda)$ and a unitary $u \in \mathcal{U}(M)$ such that $q = u z u^*$. Define $Q = u Pz u^*$. Then $Q \subset qMq$ is diffuse, amenable and $\mathcal{N}_{qMq}(Q)''$ has no amenable direct summand.

Thanks to Proposition \ref{bimodule} and since $Q \subset qMq$ is weakly compact, we may apply Theorem \ref{step} and we get $i = 1, \dots, m$ such that $Q \preceq_M M_i$. Then, there exist $n \geq 1$, a projection $p \in M_i^n$, a non-zero partial isometry $v \in \mathbf{M}_{1, n}(\C) \otimes qM$, and a unital $\ast$-homomorphism $\psi : Q \to pM_i^np$ such that $xv = v\psi(x)$, for every $x \in Q$. Note that $vv^* \in Q' \cap qMq \subset \mathcal{N}_{qMq}(Q)''$ and $v^*v \in \psi(Q)' \cap pM^np$. Note that $\psi(Q)$ is a unital von Neumann subalgebra of $pM_i^np$. We have the following alternative:

{\bf First case: assume $\psi(Q) \npreceq_{M_i^n} L(\Lambda)^n$.} If we apply Theorem $1.1$ in \cite{ipp}, we get $v^*v \in pM_i^np$ so that we may assume $v^*v = p$. Then $v^* Q v = \psi(Q)$. Take $u \in \mathcal{N}_{qMq}(Q)$. We have
\begin{eqnarray*}
(v^* u v) \psi(Q) & = & (v^* u v)v^* Q v \\
& = & v^* u Q v \\
& = & v^* Q u v \\
& = & v^*Qv(v^* u v) \\
& = & \psi(Q)(v^* u v),
\end{eqnarray*}
so that $v^* u v$ quasi-normalizes $\psi(Q)$ inside $pM^np$. Thus, since $\psi(Q) \npreceq_{M_i^n} L(\Lambda)^n$, we have $v^*\mathcal{N}_{qMq}(Q)'' v \subset pM_i^np$ by Theorem $1.1$ in \cite{ipp}. By assumption, $vv^*~\mathcal{N}_{qMq}(Q)''~vv^*$ is not amenable and 
\begin{equation*}
\Ad(v^*)(vv^*\mathcal{N}_{qMq}(Q)'' vv^*) \subset  pM_i^n p,
\end{equation*}
where $\Ad(v^*) : vv^* M vv^* \to pM^np$ is a $\ast$-isomorphism. Since $pM_i^np$ is amenable, we get a contradiction.

{\bf Second case: assume $\psi(Q) \preceq_{M_i^n} L(\Lambda)^n$.} At this point, we have $Q \preceq_M M_i$ and $\psi(Q) \preceq_{M_i^n} L(\Lambda)^n$. Our aim now is to show that $Q \preceq_{M} L(\Lambda)$. We proceed as in Remark 3.8 of \cite{vaesbimodule}. Recall that there are a projection $p \in M_i^n$, a non-zero partial isometry $v \in \mathbf{M}_{1, n}(\C) \otimes qM$, and a unital $\ast$-homomorphism $\psi : Q \to pM_i^np$ such that $xv = v\psi(x)$, for every $x \in Q$. We may assume that $p$ equals the support projection of $E_{pM_i^np}(v^*v)$. Since $\psi(Q) \preceq_{M_i^n} L(\Lambda)^n$, we get $k \geq 1$, a non-zero partial isometry $w \in p(\mathbf{M}_{n, k}(\C) \otimes M_i)$ and a (possibly non-unital) $\ast$-homomorphism $\theta : Q \to L(\Lambda)^k$ such that $\psi(x) w = w \theta(x)$, for every $x \in Q$.  Hence,
\begin{equation*}
x v w = v \psi(x) w = vw \theta(x), \forall x \in Q.
\end{equation*}
But $vw \neq 0$. Otherwise, we would have 
\begin{equation*}
E_{pM_i^np}(v^*v)w = (1 \otimes E_{pM_i^np})(v^*v w) = 0.
\end{equation*}
Since $p$ is the projection support of $E_{pM_i^np}(v^*v)$, we would get $w = pw = 0$, contradiction. Taking now the polar decomposition $vw = u |vw|$, $u$ is a non-zero partial isometry in $\mathbf{M}_{1, k}(\C) \otimes M$, such that $x u = u \theta(x)$, for any $x \in Q$. This proves that $Q \preceq_{M} L(\Lambda)$. Note that once again, $uu^* \in Q' \cap qMq \subset \mathcal{N}_{qMq}(Q)''$ and $u^*u \in \theta(Q)' \cap \theta(q)M^k\theta(q)$.

We regard now $M$ as the crossed product von Neumann algebra $L(\Upsilon) \rtimes \Lambda$. Since the $\tau$-preserving action $\Lambda \curvearrowright L(\Upsilon)$ is mixing by Proposition \ref{properties} and $\theta(Q) \subset \theta(q)L(\Lambda)^k\theta(q)$ is diffuse, it follows from Theorem $3.1$ in \cite{popamal1} (see also Theorem D.4 in \cite{vaesbern}) that $u^*u \in \theta(q)L(\Lambda)^k\theta(q)$, so that we may assume $u^*u = \theta(q)$. Note that $u^* Q u = \theta(Q)$. Moreover since $\theta(Q)$ is diffuse, Theorem 3.1 in \cite{popamal1} yields that the quasi-normalizer of $\theta(Q)$ inside $\theta(q)M^k\theta(q)$ is contained in $\theta(q) L(\Lambda)^k \theta(q)$. Proceeding exactly in the same way we did before, we get
\begin{equation*}
\Ad(u^*)(uu^*\mathcal{N}_{qMq}(Q)'' uu^*) \subset  \theta(q) L(\Lambda)^k \theta(q).
\end{equation*}
Since $\theta(q)L(\Lambda)^k \theta(q)$ is amenable and $uu^*\mathcal{N}_{qMq}(Q)'' uu^*$ is non-amenable, we finally get a contradiction, which finishes the proof.

\section{The class of groups $G$ for which $L(G)$ is strongly solid}

Denote by $\mathcal{C}_{\ssolid}$ the class of countable groups $G$ weakly amenable with contant $1$ for which the group von Neumann algebra $L(G)$ is strongly solid. This class contains all amenable groups, the free groups $\F_n$ \cite{ozawapopa}, lattices in $\SL(2, \R)$ or $\SL(2, \C)$ \cite{ozawapopaII} and the groups $\Gamma$ which satisfy Assumption A by our Theorem A. This class is obviously stable under taking subgroups. The main result of this section is that the class $\mathcal{C}_{\ssolid}$ is stable under taking free products, see Theorem B in the Introduction.

\subsection{Technical results}

We need some preparation before proving Theorem B.

\begin{lem}\label{lemmatech1}
Let $G_1, G_2$ be two non-trivial groups and write $M = L(G_1) \ast L(G_2)$ for the free product. Then there exists a diffuse abelian subalgebra $A \subset M$ such that $A \npreceq_M L(G_i)$, for any $i = 1, 2$.
\end{lem}

\begin{proof}
Let $v \in G_1 \backslash\{e\}$, $w \in G_2 \backslash\{e\}$ so that $\tau(v) = \tau(w) = 0$. Define $u = v w$ and $A$ the abelian subalgebra generated by $u$. Since $u$ is a Haar unitary, $A$ is clearly diffuse. Define $u_n = u^n$, $\forall n \in \N$. For $a, b \in G_1 \ast G_2$ be equal to $e$ or reduced words with letters alternating from $G_1$ and $G_2$, it is easy to check that for $n$ large enough $E_{L(G_i)}(a u_n b) = 0$. Using Kaplansky density theorem, we get
\begin{equation*}
\lim_{n \to \infty} \|E_{L(G_i)}(a u_n b)\|_2 = 0, \forall a, b \in M.
\end{equation*}
This means exactly that $A \npreceq_M L(G_i)$.
\end{proof}

\begin{prop}\label{proptech2}
Let $M$ be a diffuse strongly solid von Neumann algebra. Then, for any $n \geq 1$ and any non-zero projection $p \in M^n$, $pM^np$ is strongly solid.
\end{prop}

\begin{proof}
Let $M$ be a diffuse strongly solid von Neumann algebra. First assume that $n = 1$. Let $p \in M$ be a non-zero projection. Let $A \subset pMp$ be a diffuse amenable subalgebra. Let $B \subset (1 - p)M(1 - p)$ be a diffuse abelian subalgebra. Then $C = A \oplus B \subset M$ is a unital diffuse amenable subalgebra. Thus, $\mathcal{N}_M(C)''$ is amenable. Since
\begin{equation*}
\mathcal{N}_{pMp}(A)'' \oplus B \subset \mathcal{N}_M(C)'',
\end{equation*}
it follows that $\mathcal{N}_{pMp}(A)''$ is amenable, so $pMp$ is strongly solid.

Assume now that $M$ is a strongly solid ${\rm II_1}$ factor. Let $n \geq 1$ and denote by $\tau$ the canonical trace on $M^n$. Let $P \subset M^n$ be a diffuse amenable subalgebra and assume that  $\mathcal{N}_{M^n}(P)''$ is not amenable. Write $1 - z \in \mathcal{Z}(\mathcal{N}_{M^n}(P)'')$ for the maximal projection such that $\mathcal{N}_{M^n}(P)''(1 - z)$ is amenable. Then $z \neq 0$ and $\mathcal{N}_{M^n}(P)''z$ has no amenable direct summand. Since $\mathcal{N}_{M^n}(P)''z \subset \mathcal{N}_{zM^nz}(Pz)''$ is a unital inclusion (with unit $z$), $\mathcal{N}_{zM^nz}(Pz)''$ has no amenable direct summand either. If $\tau(z) \leq 1/n$, regarding $zM^nz$ as a corner of $M$ and using the first part of the proof, we obtain a contradiction because $M$ is strongly solid. If $\tau(z) \geq 1/n$, since $Pz$ is diffuse amenable we may shrink $z$ to $z_0 \in Pz$ such that $\tau(z_0) = 1/n$ and $\ctr_{Pz}(z_0) = c z_1$, where $c$ is a scalar and $z_1 \in \mathcal{Z}(Pz)$. Using Lemma $3.5$ in \cite{popamal1}, we get
\begin{equation*}
\mathcal{N}_{z_0 M^n z_0}(z_0Pz z_0)'' = z_0 \mathcal{N}_{z M^n z}(Pz)'' z_0.
\end{equation*}
We may regard $z_0 M^n z_0 \simeq M$. Let $A = z_0 Pz z_0 \subset M$. Then $A$ is a unital diffuse amenable subalgebra of $M$ and $\mathcal{N}_M(A)''$ is not amenable (because $\mathcal{N}_{z M^n z}(Pz)''$ has no amenable direct summand), which again contradicts the fact that $M$ is strongly solid.

Finally, let $M$ be a diffuse strongly solid von Neumann algebra. Let $n \geq 1$. We write the center of $M$ as 
\begin{equation*}
\mathcal{Z}(M) = \mathcal{Z}(M)z_0 \oplus \bigoplus_{i \geq 1} \C z_i,
\end{equation*}
where $\mathcal{Z}(M)z_0$ is the diffuse part. Since $Mz_0$ is strongly solid and $Mz_0 = (\mathcal{Z}(M)z_0)' \cap Mz_0$, it follows that $Mz_0$ is amenable. We obtain
\begin{equation*}
M = Mz_0 \oplus \bigoplus_{i \geq 1} Mz_i,
\end{equation*}
where $Mz_i$ are strongly solid ${\rm II_1}$ factors. Thus,
\begin{equation*}
M^n = (Mz_0)^n \oplus \bigoplus_{i \geq 1} (M z_i)^n.
\end{equation*}
Since $(Mz_0)^n$ is amenable and $(M z_i)^n$ are strongly solid ${\rm II_1}$ factors, it follows that $M^n$ is strongly solid as well, which finishes the proof.
\end{proof}

\subsection{Proof of Theorem B} 

Let $G_1, G_2 \in \mathcal{C}_{\ssolid}$. We may and will assume $|G_1| \geq 2$, $|G_2| \geq 3$. The fact that $G_1 \ast G_2$ is weakly amenable with constant $1$ follows from \cite{RicardXu}. Denote by $M = L(G_1 \ast G_2) = L(G_1) \ast L(G_2)$. Observe that since $M$ is a free product, it follows that the $M$-$M$ bimodule $L^2(\widetilde{M}) \ominus L^2(M)$ is isomorphic to a multiple of $L^2(M) \bar{\otimes} L^2(M)$. We assume that $M$ is not strongly solid and deduce a contradiction. Then there exists a diffuse amenable subalgebra $P \subset M$ such that $\mathcal{N}_M(P)''$ is not amenable. Write $1 - z \in \mathcal{Z}(\mathcal{N}_M(P)'')$ for the maximal projection such that $\mathcal{N}_M(P)''(1 - z)$ is amenable. Then $z \neq 0$ and $\mathcal{N}_M(P)''z$ has no amenable direct summand.

Let $A \subset M$ be a unital diffuse abelian subalgebra such that $A \npreceq_M L(G_i)$ as in Lemma \ref{lemmatech1}. Since $A$ is diffuse and $M$ is a ${\rm II_1}$ factor, there exists a unitary $u \in \mathcal{U}(M)$ such that $q = u z u^* \in A$. Define $Q = uPzu^*$ and $B = Q + A(1 - q)$. Observe that $B \subset M$ is a {\it unital} inclusion, $\mathcal{N}_{qMq}(Q)''$ has no amenable direct summand. Since $B \subset M$ is weakly compact and $\mathcal{N}_M(B)''$ is not amenable, applying Theorem $4.10$ in \cite{ozawapopa} or our Theorem \ref{step}, we get $i = 1, 2$ such that $B \preceq_M L(G_i)$. Thus, there exists $n \geq 1$, a non-zero partial isometry $v \in \mathbf{M}_{1, n}(\C) \otimes M$ and a (possibly non-unital) $\ast$-homomorphism $\psi : B \to L(G_i)^n$ such that $x v = v \psi(x)$, $\forall x \in B$. Observe that $q v \neq 0$, because otherwise we would have $vv^* \leq 1 - q$ and $x v  = v \psi(x)$, $\forall x \in A(1 - q)$. This would mean that $A(1 - q) \preceq_M L(G_i)$ and so $A \preceq_M L(G_i)$, contradiction. Write $q v = w |q v|$ for the polar decomposition of $qv$. It follows that $w \in \mathbf{M}_{1, n}(\C) \otimes M$ is a non-zero partial isometry such that $x w = w \psi(x)$, $\forall x \in Q$. This means exactly that $Q \preceq_M L(G_i)$. Note that $ww^* \in Q' \cap qMq \subset \mathcal{N}_{qMq}(Q)''$ and $w^*w \in \psi(Q)' \cap \psi(q) M^n \psi(q)$.

Set $N = \mathcal{N}_{qMq}(Q)''$. Since $\psi(Q)$ is a diffuse subalgebra of $\psi(q) L(G_i)^n \psi(q)$, Theorem $1.1$ in \cite{ipp} yields $w^* N w \subset \psi(q) L(G_i)^n \psi(q)$. We may and will assume that $w^*w = \psi(q)$. Let $\widetilde{z} \in \mathcal{Z}(N)$ be a projection such that $\widetilde{z} = \sum_{j = 1}^k v_j v_j^*$ with $v_j$ partial isometries in $N$ and $v_j^* v_j \leq ww^*$. Define the non-zero partial isometry
\begin{equation*}
\widetilde{w} = [v_1 w \cdots v_k w] \in \mathbf{M}_{1, kn}(\C) \otimes M.
\end{equation*}
Observe that $\widetilde{z} = \widetilde{w}\widetilde{w}^* \in \mathcal{Z}(N)$ and $\widetilde{p} = \widetilde{w}^*\widetilde{w} \in L(G_i)^{kn}$. Then $\widetilde{Q} = \widetilde{w}^* Q \widetilde{w}$ is a unital diffuse amenable subalgebra of $\widetilde{p} L(G_i)^{kn} \widetilde{p}$ and for any $u \in \mathcal{N}_{qMq}(Q)$, $\widetilde{w}^* u \widetilde{w} \in \mathcal{N}_{\widetilde{p}L(G_i)^{kn}\widetilde{p}}(\widetilde{Q})$, because $\widetilde{z} = \widetilde{w} \widetilde{w}^* \in \mathcal{Z}(N)$. We obtain
\begin{equation*}
\Ad(\widetilde{w}^*)(N \widetilde{z}) \subset \mathcal{N}_{\widetilde{p}L(G_i)^{kn}\widetilde{p}}(\widetilde{Q})'',
\end{equation*}
where $\Ad(\widetilde{w}^*) : \widetilde{z} M \widetilde{z} \to \widetilde{p} M^{kn} \widetilde{p}$ is a $\ast$-isomorphism. Since $N \widetilde{z}$ is not amenable, it follows that $\mathcal{N}_{\widetilde{p}L(G_i)^{kn}\widetilde{p}}(\widetilde{Q})''$ is not amenable either. This contradicts the fact that $\widetilde{p}L(G_i)^{kn}\widetilde{p}$ is strongly solid (see Proposition \ref{proptech2}), which finishes the proof.

\section{${\rm II_1}$ factors with at most one Cartan subalgebra}

\subsection{Crossed products as amalgamated free products} 

Let $G_1, G_2$ be countable groups and denote by $G = G_1 \ast G_2$ their free product. Let $(Q, \tau)$ be an amenable finite von Neumann algebra and $G \curvearrowright Q$ be a $\tau$-preserving action. Denote by $M = Q \rtimes G$ the corresponding crossed product von Neumann algebra and observe that $M$ may be regarded as the amalgamated free product
\begin{equation*}
M = (Q \rtimes G_1) \ast_Q (Q \rtimes G_2).
\end{equation*}
Using the notation of Section \ref{keyresult}, we have
\begin{equation*}
\widetilde{M} = Q \rtimes (G \ast \F_2),
\end{equation*}
where $\F_2$ acts trivially on $Q$. The next proposition is an easy consequence of Lemma $4.7$ in \cite{ozawapopa}:

\begin{prop}\label{decompositionbis}
The $M$-$M$ bimodule $L^2(\widetilde{M}) \ominus L^2(M)$ is isomorphic to a multiple of $L^2\langle M, e_Q\rangle$.
\end{prop}

\subsection{Proof of Theorem C} 

Let $G_1, G_2$ be weakly amenable groups with constant $1$, so that $G = G_1 \ast G_2$ is still weakly amenable with constant $1$ (see \cite{RicardXu}). We assume $|G_1| \geq 2$, $|G_2| \geq 3$, so that $G$ is not  amenable. Let $(Q, \tau)$ be a finite amenable von Neumann algebra together with a $\tau$-preserving action $G \curvearrowright Q$. We will assume: 
\begin{enumerate}
\item Either $Q = \C$.
\item Or $Q = L^\infty(X, \mu)$ and the p.m.p. action $G \curvearrowright (X, \mu)$ is assumed to be free ergodic and {\it profinite}. Note that this assumption forces $G_1, G_2$ to be residually finite.
\end{enumerate}
We shall denote by $M = Q \rtimes G$ the crossed product ${\rm II_1}$ factor which can be regarded as
\begin{equation*}
M = M_1 \ast_Q M_2
\end{equation*}
where $M_i = Q \rtimes G_i$. Note that under the assumptions $(1)$ or $(2)$, $M$ always has the c.m.a.p., i.e. $\Lambda_{\cb}(M) = 1$. 

Let $A \subset M$ be a Cartan subalgebra, i.e. $A = A' \cap M$ and $\mathcal{N}_M(A)'' = M$. Thanks to Proposition \ref{decompositionbis} and since $A \subset M$ is weakly compact, Theorem \ref{step} yields $i = 1, 2$ such that $A \preceq_M M_i$. Thus there exist $n \geq 1$, a non-zero partial isometry $v \in \mathbf{M}_{1, n}(\C) \otimes M$, a projection $p \in M_i^n$ and a unital $\ast$-homomorphism $\psi : A \to pM_i^np$ such that $x v = v \psi(x)$, $\forall x \in A$. We may and will assume that the support projection of $E_{M_i^n}(v^*v)$ equals $p$. As in the proof of Theorem A, we have the following alternative.

{\bf Assume that $\psi(A) \npreceq_{M_i^n} Q^n$.} If we apply Theorem $1.1$ in \cite{ipp}, we may assume $p = v^*v$ and we have
\begin{equation*}
v^* M v = v^* \mathcal{N}_M(A)'' v \subset pM_i^np.
\end{equation*}
Thus, we get $pM^np = pM_i^np$. Since $M^n = M_1^n \ast_{Q^n} M_2^n$ (see Proposition \ref{decompositionbis}), we have
\begin{equation*}
L^2(M^n) \ominus L^2(M_i^n) \cong \bigoplus L^2(M_i^n) \bar{\otimes}_{Q^n} L^2(M_i^n).
\end{equation*} 
as $M_i^n$-$M_i^n$ bimodules. Consequently, we have $L^2(p M_i^n) \bar{\otimes}_{Q^n} L^2(M_i^n p) = p(L^2(M_i^n) \bar{\otimes}_{Q^n} L^2(M_i^n))p = 0$, contradiction.

{\bf Assume that $\psi(A) \preceq_{M_i^n} Q^n$.} At this point we have $A \preceq_M M_i$ and $\psi(A) \preceq_{M_i^n} Q^n$. As in the proof of Theorem A, we obtain that $A \preceq_M Q$. 
\begin{enumerate}
\item If $Q = \C$, we get a contradiction since $A$ is diffuse. This yields that $L(G)$ has no Cartan subalgebra.
\item If $Q = L^\infty(X, \mu)$ is the Cartan subalgebra coming from the profinite action of $G$ on $(X, \mu)$, then applying Theorem A.1 in \cite{popa2001}, we obtain $u \in \mathcal{U}(M)$ such that $u A u^* = L^\infty(X, \mu)$. This yields that $L^\infty(X, \mu) \subset L^\infty(X, \mu) \rtimes G$ is the unique Cartan subalgebra up to unitary conjugacy.
\end{enumerate}

%\addcontentsline{toc}{section}{Bibliography}
\bibliographystyle{plain}

\end{document}